\documentclass[10pt]{amsart}

\usepackage{amsmath, amsthm,amscd, amssymb}

\usepackage{MnSymbol}

\usepackage[dvips]{graphicx}
\usepackage{color}
\usepackage[all]{xy}
\usepackage{verbatim}

\newtheorem{Lemma1}{{Lemma}}[section]
\newtheorem{Theo1}[Lemma1]{{Theorem}}
\newtheorem*{Theo2}{{Theorem}}
\newtheorem{Def1}[Lemma1]{{Definition}}
\newtheorem{Prop1}[Lemma1]{{Proposition}}
\newtheorem{Claim1}[Lemma1]{{Claim}}
\newtheorem{Rem1}[Lemma1]{{Remark}}
\newtheorem{Cor1}[Lemma1]{{Corollary}}
\newtheorem{Ex1}[Lemma1]{{Example}}
\newtheorem{Not1}[Lemma1]{{Notation}}

\newenvironment{Lemma}{\begin{Lemma1}}{\end{Lemma1}}
\newenvironment{Def}{\begin{Def1}\rm}{\end{Def1}}
\newenvironment{Prop}{\begin{Prop1}}{\end{Prop1}}
\newenvironment{Rem}{\begin{Rem1}\rm}{\end{Rem1}}
\newenvironment{Theorem}{\begin{Theo1}}{\end{Theo1}}

\newenvironment{Cor}{\begin{Cor1}}{\end{Cor1}}

\newenvironment{Example}{\begin{Ex1}\rm}{\end{Ex1}}

\title{Dg-separable dg-extensions}

\author{Alexander Zimmermann}
\address{\newline
Universit\'e de Picardie,
\newline D\'epartement de Math\'ematiques et LAMFA (UMR 7352 du CNRS),
\newline 33 rue St Leu,
\newline F-80039 Amiens Cedex 1,
\newline France}
\email{alexander.zimmermann@u-picardie.fr\newline  {\tt https://orcid.org/orcid/0000-0001-8168-4624}}

\date{December 2, 2024; revised December 26, 2025}

\newcommand{\lra}{\longrightarrow}

\newcommand{\ra}{\rightarrow}
\newcommand{\sdp}{\times\kern-.2em\vrule height1.1ex depth-.05ex}
\newcommand{\epi}{\lra \kern-.8em\ra}

\newcommand{\Z}{{\mathbb Z}}

\newcommand{\im}{\textup{im}}
\newcommand{\dickebox}{{\vrule height5pt width5pt depth0pt}}

\newcommand{\id}{\textup{id}}

\newcommand{\Hom}{\textup{Hom}}

\newcommand{\End}{\textup{End}}

 \normalbaselines
\subjclass[2020]{Primary: 16E45; Secondary: 12H05; 13A02;   16W50}
\keywords{differential graded algebras}

\newcommand{\red}[1]{}

\begin{document}

\begin{abstract}
We define and characterise completely dg-separable dg-extensions
$\varphi:(A,d_A)\ra(B,d_B)$.
We completely characterise the case of graded commutative dg-division algebras
in characteristic different from $2$. 
We prove that for a dg-separable extension a short exact sequence 
of dg-modules over $(B,d_B)$ splits if and only if the restriction 
to $(A,d_A)$ splits.   
\end{abstract}

\maketitle

\section*{Introduction}

Let $K$ be a commutative ring.
A differential graded $K$-algebra $(A,d)$ is a $\Z$-graded algebra $A$ together with a
$K$-linear graded endomorphism $d:A\ra A$ of degree $1$ such that $d^2=0$ and
$$d(a\cdot b)=d(a)\cdot b+(-1)^{|a|}a\cdot d(b)$$ 
for all homogeneous $a,b\in A$, where we denote by $|a|$ the degree of $a\in A$.
Differential graded algebras (or dg-algebras for short) 
were defined by Cartan~\cite{Cartandg} in 1954
and proved to be highly successful in many subjects, such as homological algebra,
algebraic topology, differential and algebraic geometry, and alike.
However, the ring theory of differential graded algebras remained largely unexplored 
until quite recently. The first of the results in this direction was a characterisation
of Aldrich and Garcia Rozas \cite{Tempest-Garcia-Rochas} of acyclic dg-algebras.
\cite{dgorders} then studied general ring theoretic properties, such as a dg-Nakayama lemma, and
independently in a parallel development 
Orlov \cite{Orlov1} studied finite dimensional dg-$K$-algebras over a field $K$. Goodbody 
\cite{Goodbody} proved a version of Nakayama's lemma in the dg-setting following Orlov's approach. 
In the sequel \cite{dgBrauer} defined and studied a dg-Brauer group, and in \cite{dgGoldie} 
Ore localisation and a Goldie theorem was studied in the context of dg-algebras. Further, 
in \cite{dgfields} a concept of a dg-division algebra was developed, and a complete classification 
was given.  In this case we showed that a dg-division algebra is 
either acyclic or has differential $d=0$. Note that in \cite{dgBrauer} a technical hypothesis 
was imposed for the classification. We shall prove in this paper that the technical assumption is superfluous.

In general, a $K$-algebra $B$ is called separable over a $K$-subalgebra $A$
if the multiplication map $B\otimes_AB\lra B$ is split as morphism of $B-B$-bimodules. 
A graded version was given by N\u ast\u asescu-van Oystaeyen~\cite{gradedrings}
asking for a split in the category of graded bimodules. 

We define in this paper a differential graded separability, asking simply that the splitting 
of the multiplication map is a map of differential graded bimodules. 
We use the classification from   \cite{dgfields} to show that a field extension between two 
graded-commutative acyclic dg-division rings is dg-separable if the extension of cycles is 
 graded-separable. In characteristic different from $2$, the converse also holds. 
Further, we show that in characteristic different from $2$ an extension of dg-division algebras 
from an algebra with differential $0$ to an acyclic algebra is never separable. We finally note that
a dg-extension $(A,d_A)\ra (B,d_B)$ where $(A,d_A)$ is acyclic implies that $(B,d_B)$ is acyclic as well. 
This gives a complete picture of separability of extensions of graded commutative dg-division algebras
in characteristic different from $2$. The results are displayed in Theorem~\ref{maintheorem}. 

In general, we show that a dg-extension $\varphi:(A,d_A)\ra (B,d_B)$ 
is dg-separable if and only if there is a homogeneous element 
$\omega\in\ker(d_{B\otimes_AB})$ of degree $0$ 
with $b\omega=\omega b$ for all $b\in B$ and 
mapping to $1$ under the multiplication map $B\otimes_AB\ra B$. 
We show in Theorem~\ref{relativelysplit} 
that this then implies that a short exact sequence of dg-modules over $(B,d_B)$
if and only if the restriction to $(A,d_A)$ splits.  We further mention that our 
concept of dg-separability gives that the restriction functor is a separable functor
in the sense of N\u ast\u asescu, van den Bergh and van Oystaeyen~\cite{NastacescuvdBerghvOystaen}.

The paper is organised as follows. In Section~\ref{dgdivisioalgebrasrevis} we recall 
results from \cite{dgfields} concerning dg-division rings as far as they are relevant for 
this work. Section~\ref{dgseparbilitysect} then gives the definition of a 
dg-separable extension. In Section~\ref{dgseparablefieldextension} we completely classify
dg-separable dg-extension of graded-commutative dg-division rings, which includes our first
main result Theorem~\ref{maintheorem}. 
Finally, Section~\ref{generalpropertiesofdgseparableext} shows the second main  
result Theorem~\ref{relativelysplit}.

\noindent
{\bf Acknowledgement:} I would like to thank Xiao-Wu Chen for his kind invitation to 
USTC. This opportunity enabled me to improve the paper. I would also like to thank the referees for their time and effort in reviewing the manuscript.

\section{Dg-division algebras revisited}

\label{dgdivisioalgebrasrevis}

First recall some notations. As a reference one may take \cite{Yekutielibook} or \cite{dgorders,dgGoldie,dgfields}.
Let $(A,d)$ be a dg-$K$-algebra. Then a left dg-module over $(A,d)$ is a $\Z$-graded $A$-module $M$
together with an endomorphism $\delta:M\ra M$ of degree $1$ with $\delta^2=0$ and 
$\delta(a\cdot m)=d(a)\cdot m+(-1)^{|a|}a\cdot\delta(m)$ for all homogeneous $a\in A$ and $m\in M$. 
If $(A,d)$ is a dg-algebra, then $(A^{op},d)$ is a dg-algebra as well, where $A^{op}$ coincides with
$A$ as $K$-module, and where $a\cdot_{op}b:=(-1)^{|a||b|}b\cdot a$ for all homogeneous $a,b\in A$. 
Further, a right dg-module over $(A,d)$ is a left dg-module over $(A^{op},d)$. 
For two dg-modules $(M,\delta_M)$ and $(N,\delta_N)$ over $(A,d)$ we set 
$$\Hom_A^k((M,\delta_M),(N,\delta_N)):=\{f\in\Hom_{K,\textup{ graded}}(M,N)\;|\;f(am)=(-1)^{|a|k}af(m)\}$$
and put $d_\Hom(f):=\delta_M\circ f-(-1)^{|f|}\circ\delta_M$.

We abbreviate 
$$\Hom^\bullet_A((M,\delta_M),(N,\delta_N)):=\bigoplus_{k\in\Z}\Hom^k_A((M,\delta_M),(N,\delta_N)).$$
and 
$$\End^\bullet_A((M,\delta_M)):=\Hom_A^\bullet((M,\delta_M),(M,\delta_M)). $$

Then $\End^\bullet_A( (M,\delta_M)),d_\Hom)$ is a dg-algebra, and 
$\Hom^\bullet_A( (M,\delta_M),(N,\delta_N)),d_\Hom)$ is a dg-bimodule over
$\End^\bullet_A( (N,\delta_N)),d_\Hom)$- $\End^\bullet_A( (M,\delta_M)),d_\Hom)$.

Recall from \cite{dgfields} the definition of a differential graded division algebra.

\begin{Def} \cite{dgfields}
A {\em dg-division algebra} is a dg-algebra $(A,d)$ such that the only dg-left ideals are $0$ and $A$ and
the only dg-right ideals are $0$ and $A$.  
\end{Def}

Differential graded division algebras $(A,d)$ 
were completely classified in \cite{dgfields} (cf Remark~\ref{notnecessaryrem} below).

The corresponding result is the following. 

\begin{Theorem} \cite{dgfields} \label{dgdivisionalgebraclassification}
Let $(A,d)$ be a dg-algebra. 
Then 
\begin{itemize}
\item
$(A,d)$ is a dg-division algebra if and only if $\ker(d)$ is a 
$\Z$-gr-division algebra (cf \cite{gradedrings}).
\item 
If $(A,d)$ is a dg-division algebra,
\begin{itemize} 
\item then,
\begin{itemize}
\item  either $d=0$ and $\ker(d)$ is  a skew-field 
concentrated in degree $0$, 
\item or $H(A,d)=0$ and 
\begin{itemize}
\item either $\ker(d)$ is a skew-field concentrated in degree $0$
\item or there is a skew field $R_0$
such that $\ker(d)\simeq R_0[X,X^{-1};\phi]$ for an automorphism $\phi$ of $R_0$ 
and $Xr=\phi(r)X$ for any $r\in R_0$. 
\end{itemize}
\end{itemize}
\item If $H(A,d)=0$, then there is a homogeneous element $y$ with
$d(y)=1$ and $y^2\in\ker(d)$, and there is a map $D:\ker(d)\lra\ker(d)$  of degree $1$
defined by  $$D(a)=-(-1)^{|a|}d(yay)=ya-(-1)^{|a|}ay$$ 
for any homogeneous $a\in\ker(d)$,  
such that  $A$ is isomorphic with the quotient of the twisted polynomial ring
$$A\simeq\ker(d)[T;D]/(T^2-y^2).$$ 
Moreover, the algebra structure on the twisted group ring is 
given by $D(a)=Ta-(-1)^{|a|}aT$ for any homogeneous $a\in\ker(d)$. Furthermore, 
$A=\ker(d)\oplus y\ker(d)$, and the isomorphism 
is 
\begin{eqnarray*}
\Phi:\ker(d)[T;D]&\lra& A\\
b+Ta&\mapsto&b+ya
\end{eqnarray*} 
for any homogeneous $a,b\in\ker(d)$. Further, 
for any homogeneous $a,b\in\ker(d)$ we get $d(b+ya)=a$. 
\end{itemize}
\end{itemize}
\end{Theorem}

\begin{Rem}  \label{notnecessaryrem}
Recall that a $\Z$-graded ring is a $\Z$-graded-division ring if every homogeneous element is invertible. 
In \cite{dgfields} an additional hypothesis on $(A,d)$ was imposed in the formulation of Theorem~\ref{dgdivisionalgebraclassification}. Namely we asked 
that  the set of left regular homogeneous 
elements of $\ker(d)$ coincides with the set of right regular homogeneous elements of $\ker(d)$.
This was used in order to show that $\ker(d)$ is a graded-division ring. 
We shall show here that this hypothesis is unnecessary.
\end{Rem}

\begin{Lemma}\label{dgdivisionalgequalsgrdivalg}
A dg-algebra $A$ is a dg-division ring if and only if $\ker(d)$ is a $\Z$-graded-division ring.
\end{Lemma}

Proof. If $\ker(d)$ is a $\Z$-graded-division ring, then $(A,d)$ is a dg-division algebra by 
\cite[Lemma 2.1]{dgfields}. 

Suppose that $\ker(d)$ is not a $\Z$-graded-division ring, and 
that $(A,d)$ is a dg-division ring. Then there is 
an in $\ker(d)$ non invertible homogeneous element $0\neq x\in\ker(d)$. 
Hence, $x$ is not left invertible, or not right 
invertible in $\ker(d)$, or both, since else 
$x$ would be invertible, contradicting the hypothesis. 

Suppose first that $x\cdot\ker(d)\neq\ker(d)$.
If $xA\neq A$, then, since $x\in\ker(d)$, the ideal 
$xA$ would be a non trivial dg-right ideal of $A$. This is impossible since 
$(A,d)$ is assumed to be a dg-division algebra.
Hence $xA=A$, and therefore there is a homogeneous $y\in A$ with $xy=1$.
But then 
$$0=d(1)=d(xy)=d(x)\cdot y+(-1)^{|x|}x\cdot d(y)=(-1)^{|x|}x\cdot d(y).$$
Since $x\in\ker(d)$, the left ideal $Ax$ is a dg-left ideal of $A$. 
Hence, since $(A,d)$ is a dg-division ring, $Ax=A$. But then 
$$A\cdot d(y)=A\cdot x\cdot d(y)=A\cdot 0 =0,$$
which implies $d(y)=0$. Hence $y\in\ker(d)$, which implies in turn 
$$\ker(d)\supseteq x\cdot\ker(d)\supseteq x\cdot y\cdot \ker(d)=\ker(d).$$ 
This was excluded.

If $\ker(d)\cdot x\neq \ker(d)$, the analogous argument gives a contradiction as well. Hence $\ker(d)$ is a $\Z-gr$-division ring. 
\dickebox

\begin{Cor}
Let $(A,d)$ be a dg-algebra. 
Suppose that $(A,d)$ is a graded commutative dg-division algebra. 
Then either $\ker(d)$ is a field concentrated in degree $0$, or else 
$\ker(d)=K[X,X^{-1}]$ for some field $K$ and $X$ in non zero degree, 
and, in case $K$ is of characteristic different from $2$, then $X$ is in even degree.
\end{Cor}

We shall need to recall the definition of a differential graded structure on a
tensor product of algebras. Let $(A,d_A)$ be a dg-algebra and let $(B,d_B)$ 
de a dg-algebra. Consider a dg-homomorphism $(A,d_A)\lra (B,d_B)$. Then 
for $$d_{B\otimes_AB}=d_B\otimes \textup{id}_B+\textup{id}_B\otimes d_B,$$
respecting the Koszul sign rule,
defines a dg-$B-B$-bimodule structure on $(B\otimes_AB)$. If $A$ is a subalgebra of the
graded centre of $B$, then $(B\otimes_AB,d_{B\otimes_AB})$ is a dg-algebra again. 

\section{Dg-separability}

\label{dgseparbilitysect}

Recall that an algebra $A$ is separable if $A$ is a projective object 
in the category of $A-A$-bimodules. This is equivalent with the fact that the 
multiplication map is split as morphism of $A-A$-bimodules. 
Similarly, a graded algebra $A$ is graded separable if 
the graded bimodule $A$ is projective in the category of graded bimodules. 

\begin{Prop} \cite[Example 2.5]{ElFadilQuebec}
The extension of graded rings 
$R[T^n,T^{-n}]\subseteq S[T,T^{-1}]$  is graded-separable if and only if 
the extension $R\subseteq S$ is separable and $n$ is invertible in $R$.
\end{Prop}

Recall (cf e.g. \cite{reptheobuch}) that an algebra extension $\beta:A\lra B$ 
of $K$-algebras is separable if the multiplication map
$$\mu:B\otimes_AB\lra B$$
splits as a homomorphism of $B\otimes_K B$-bimodules. 

We shall use the analogous concept. 

\begin{Def}
Let $K$ be a commutative ring and let $(A,d_A)$ and $(B,d_B)$ be 
differential graded algebras. A {\em dg-extension of dg-algebras} is a homomorphism
$\beta:(A,d_A)\lra (B,d_B)$ 
dg-algebras. 

An extension of dg-algebras  $\beta:(A,d_A)\lra (B,d_B)$  is called {\em dg-separable} if  
the multiplication map 
$$\mu:(B,d_B)\otimes_A(B,d_B)\lra (B,d_B)$$ is split as 
morphism of differential graded $B-B$-bimodules. 
\end{Def}

Note however, for
$$B\otimes_AB\stackrel\mu\lra B$$
one has $d_{B\otimes_AB}=d_B\otimes 1+1\otimes d_B$ and then, by the Leibniz formula and the Koszul 
sign rule on graded rings, we always have 
$$\mu\circ d_{B\otimes_AB}=d_B\circ \mu.$$

\begin{Prop}\label{dgcasimir}
Let $(A,d_A)\lra (B,d_B)$ be a dg-extension of dg-algebras. This extension is dg-separable
if and only if there is $\omega\in \ker(d_{B\otimes_AB})$ homogeneous of 
degree $0$ with $b\omega=\omega b$
for all $b\in B$ and  $\mu(\omega)=1$. 
\end{Prop}

Proof.
Let $\rho:B\lra B\otimes_AB$ be a retract with $\mu\circ\rho=1_B$.  
Then $$\rho\circ d_B=d_{B\otimes_AB}\circ\rho$$
is equivalent with 
$$d_{B\otimes_AB}(\omega)=0$$
for $b\omega=\omega b$ and $\mu(\omega)=1$. 
Hence, $\omega$ has to be a cycle in $B\otimes_AB$. Even better, this is equivalent.
Suppose $\omega\in\ker(d_{B\otimes_AB})$ with $b\omega=\omega b$ for all $b$ and $\mu(\omega)=1$.
Then 
\begin{eqnarray*}
d_{B\otimes_AB}(\rho(b))&=&d_{B\otimes_AB}(b\rho(1))\\
&=&d_{B\otimes_AB}(b\omega)\\
&=&((d_B\otimes 1)(b\otimes 1)+(1\otimes d_B)(b\otimes 1))\omega+(-1)^{|b|}b d_{B\otimes B}(\omega))\\
&=&d_B(b)\omega\\
&=&d_B(b)\rho(1)\\
&=&\rho(d_B(b))
\end{eqnarray*}
where the last equation holds since $\rho$ is a morphism of bimodules. 
\dickebox

\begin{Lemma}\label{extensionoddgalgsinduceextofcycles}
Let $(A,d_A)$ and $(B,d_B)$ be dg-algebras and let $\varphi:(A,d_A)\lra (B,d_B)$ be 
a dg-extension of dg-algebras. Then $\varphi|_{\ker(d_A)}$ is an extension of 
graded rings $\ker(d_A)\lra\ker(d_B)$. 
\end{Lemma}

Proof. 
Suppose that $\varphi:(A,d_A)\lra (B,d_B)$ is a dg-extension of dg-algebras. 
Then $\varphi$ induces a graded-extension $\varphi|_{\ker(d_A)}:\ker(d_A)\lra\ker(d_B)$  
by restriction. Indeed, if $d_A(x)=0$, then 
$$0=\varphi(d_A(x))=d_B(\varphi(x))$$
and hence $\varphi(x)\in\ker(d_B)$ as well. \dickebox

\section{Characterisation of dg-separable dg-field extensions}

\label{dgseparablefieldextension}

\begin{Prop}\label{gradedseparableofcyclesimpliesseparable}
Suppose that $(A,d_A)$ and $(B,d_B)$ are graded commutative
dg-division algebras, suppose that $(A,d_A)$ is acyclic, 
and suppose that $\varphi:(A,d_A)\lra (B,d_B)$ is a dg-extension of dg-algebras.
If the restriction  $\varphi|_{\ker(d_A)}$ is a graded-separable extension, then
$\varphi$ is a dg-separable extension. If the characteristic of $A$ is different from $2$, 
then the converse also holds. 
\end{Prop}

Proof.  
The algebra $(B,d_B)$ is a left dg-module over $(A,d_A)$ via $\varphi$.  Therefore,
 by \cite{Tempest-Garcia-Rochas}
we get that $(B,d_B)$ is acyclic as well. We may hence suppose that 
$(A,d_A)$ and $(B,d_B)$ are both acyclic dg-division algebras.
Then 
$$A=\ker(d_A)[T;D_A]/(T^2-y_A^2)$$
and 
$$B=\ker(d_B)[T,D_B]/(T^2-y_B^2)$$
for $d_A(y_A)=1$ and $d_B(y_B)=1$. Suppose that $\varphi:(A,d_A)\lra (B,d_B)$ is a 
dg-homomorphism. 
Further, 
$\varphi(1_A)=1_B$ implies that we may assume that $\varphi(y_A)=y_B$. 
 
By Lemma~\ref{extensionoddgalgsinduceextofcycles} the restriction of $\varphi$ to $\ker(d_A)$ is an 
extension of graded rings $\ker(d_A)\lra\ker(d_B)$. 
 
\medskip 
 
Suppose now that $\varphi|_{\ker(d_A)}$ is a graded-separable extension. 
Let $\omega_{\ker}$ be the element from Proposition~\ref{dgcasimir} with $\mu(\omega_{\ker})=1$
and $b\cdot\omega_{\ker}=\omega_{\ker}\cdot b$ for all homogeneous $b\in\ker(d_B)$. 
Recall $$A=\ker(d_A)[T;D_A]/(T^2-y_A^2).$$
But, $\varphi(T)$ can be used as $T$ in the isomorphism 
$$B=\ker(d_B)[T,D_B]/(T^2-y_B^2)$$
since we can put $\varphi(y_A)=y_B$. 
But than we only need to show 
$$T\omega_{\ker}=\omega_{\ker} T.$$
If $(B,d_B)$ is graded commutative, then for all homogeneous $b_1,b_2\in\ker(d_B)$ 
we have that either the characteristic is $2$, or else all $b_i$ are of even degree. Hence
$$T(b_1\otimes_Ab_2)=b_1T\otimes_Ab_2=b_1\otimes_ATb_2=b_1\otimes_Ab_2T=(b_1\otimes_Ab_2)T.$$
Moreover, trivially, $T$ commutes with $T$, and hence $T$ commutes with any 
$(b_1\otimes_Ab_2)\in B\otimes_AB$. 

\medskip

As for the converse, suppose that $(A,d_A)\lra (B,d_B)$ be a dg-separable dg-extension. 
Therefore, by Proposition~\ref{dgcasimir} there is an element $\omega\in B\otimes_AB$ of degree $0$. 
with $b\omega=\omega b$ for all homogeneous $b\in B$, and the image of $\omega$ under 
the multiplication map $B\otimes_AB\ra B$ is $1$. Further, $\omega\in\ker(d_{B\otimes_AB})$. 
If $B$ is a graded commutative dg-division ring of characteristic different from $2$, 
then $\ker(d_B)$ has to be concentrated
in even degrees, since any homogeneous element is invertible, whence not nilpotent ($x^2=-x^2$ for
elements of odd degree),  and by consequence in even degrees. 

Consider the map 
$$\Upsilon:\ker(d_B)\otimes_{\ker(d_A)}\ker(d_B)\lra B\otimes_AB$$
given by the natural inclusion. 
But then, as $B=\ker(d_B)\oplus T\ker(d_B)=\ker(d_B)\oplus\ker(d_B)T$,
and since $T$ is of degree $-1$, we see that the direct summand 
$\ker(d_B)\otimes_A T\ker(d_B)$ and $T\ker(d_B)\otimes_A\ker(d_B)$ are in odd degrees. 
Hence the image of $\Upsilon$ is in the subspace of even degree of $B\otimes_AB$. 
Further, all $b_1\otimes b_2$ with $b_1,b_2\in\ker(d_B)$ are in the image of $\Upsilon.$
Also, for all $b_1,b_2\in\ker(d_B)$ the elements 
$$Tb_1\otimes Tb_2=T^2b_1\otimes b_2=y_B^2b_1\otimes b_2$$
are in the image of $\Upsilon$. Further, $y_B^2$ is homogeneous of degree $-2$, satisfying $d_B(y_B)=1$, 
and since $$d_B(y_B^2)=d_B(y_B)y_B-y_B d_B(y_B)=y_B-y_B=0$$ 
we have in $y_B^2\in\ker(d_B)$. 
We have two cases. If $y_B^2=0$, then 
$$T\ker(d_B)\otimes_AT\ker(d_B)=y_B^2\ker(d_B)\otimes_A\ker(d_B)=0.$$
Else, $y_B^2\in\ker(d_B)^\times$ since $(B,d_B)$ is an acyclic 
dg-division algebra (cf Lemma~\ref{dgdivisionalgequalsgrdivalg}).  
Hence, the image of $\Upsilon$ is precisely the 
subspace of even degree elements of $B\otimes_AB$. 

Since $\omega$ has to be homogeneous 
of degree $0$, which is even, $\omega\in\im(\Upsilon)$. Let $\omega'\in\ker(d_B)\otimes_A\ker(d_B)$
with $\Upsilon(\omega')=\omega$. Note that $\Upsilon$ is injective. 
Therefore $\omega'\in\ker(d_B)\otimes_A\ker(d_B)$ can be used 
as required element to show that $\ker(d_A)\lra\ker(d_B)$ is graded separable. 
\dickebox



\begin{Rem}
The case of differential $0$ is trivial. A dg-extension  $(A,0)\lra (B,0)$ is precisely a graded extension. 
Note that if $(A,d_A)$ is acyclic, then any extension $(B,d_B)$ of $(A,d_A)$ is acylic as well. 
The only case left is when $d_A=0$ and $(B,d_B)$ is acyclic. 
\end{Rem}

\begin{Prop} \label{cyclestoacyclic}
Let $(A,0)$ be a gr-field and let $(B,d)$ be a graded commutative acyclic dg-division algebra. 
Suppose that $A\lra B$ is a dg-extension. If $B$ is of characteristic different from $2$, then 
this extension is not dg-separable. 
\end{Prop}

Proof. We need to find an element $\omega\in B\otimes_AB$, homogeneous of degree $0$ and 
mapping to $1$ under the multiplication map $B\otimes_AB\lra B$, such that $d_{B\otimes_BA}(\omega)=0$
and  such that $b\omega=\omega b$ for all homogeneous  $b\in B$.

If $B$ is a graded commutative dg-division ring over $k$, then $\ker(d)$ has to be concentrated
in even degrees, since any homogeneous element is invertible, whence not nilpotent ($x^2=-x^2$ for
elements of odd degree), 
and by consequence in even degrees. But then, as $B=\ker(d)\oplus T\ker(d)=\ker(d)\oplus\ker(d)T$,
and since $T$ is of degree $-1$, we see that for the element $\omega$ we have 
$$\omega\in(\ker(d)\otimes\ker(d))\oplus (T\ker(d)\otimes T\ker(d)).$$  
However, $\omega\in\ker(d_{B\otimes_AB})$. Clearly, 
$$(\ker(d)\otimes\ker(d))\subseteq \ker(d_{B\otimes_AB}).$$
Now, for $x=\sum_{i=1}^nTb_i\otimes Tb_i'\in (T\ker(d)\otimes T\ker(d))$
we get 
$$
0=d_{B\otimes_AB}(x)=\sum_{i=1}^nb_i\otimes Tb_i'+(-1)^{|b_i|+1}Tb_i\otimes b_i'
=\sum_{i=1}^nb_i\otimes Tb_i'-Tb_i\otimes b_i'
$$
since all $b_i$ are of even degree. 
Since, 
$$B\otimes_AB=(\ker(d)\otimes\ker(d))\oplus(T\ker(d)\otimes\ker(d))
\oplus(\ker(d)\otimes T\ker(d))\oplus (T\ker(d)\otimes T\ker(d))$$
we get $$(T\ker(d)\otimes\ker(d))\cap(\ker(d)\otimes T\ker(d))=0$$
and hence 
$$\sum_{i=1}^nb_i\otimes Tb_i'=0=\sum_{i=1}^nTb_i\otimes b_i'=T\cdot(\sum_{i=1}^nb_i\otimes b_i')$$
which shows that $$\sum_{i=1}^nb_i\otimes b_i'=0.$$
Therefore $x=0$.
But for $\omega\in\ker(d)\otimes_A\ker(d)$ we get that $T\omega\neq\omega T$ since the left hand side is 
in $T\ker(d)\otimes_A\ker(d)$ and the right hand side lies in $\ker(d)\otimes_A\ker(d)T$, whose intersection is $0$. 
\dickebox

\medskip

We illustrate the argument by a simple  

\begin{Example}\begin{itemize}\item
We illustrate the proof of Proposition~\ref{cyclestoacyclic} with an example. 
Let $A=K[X]/X^2$ for $d(X)=1$, and $X$ in degree $-1$. 
Here, $\ker(d)=K\cdot 1$.
Then we need to see if the multiplication map
\begin{eqnarray*}
A\otimes_KA&\lra& A\\
(a+ bX)\otimes(c+dX)&\mapsto&ac+(ad+bc)X
\end{eqnarray*} 
is split as $A-A$-dg-bimodules. 
As we have a $K$-basis $\{1,X\}$ of $A$, we also have a $K$-basis $\{1\otimes 1,1\otimes X, X\otimes 1,X\otimes X\}$
of $A\otimes_KA$.
The multiplication map splits as dg-map if and only if there is an element $\omega$ of degree
$0$ with $v\omega=\omega v$ for all $v\in A$ and mapping to $1$ under the multiplication map. 
As $\ker(d)$ is commutative and central, we only need to 
verify this property for $v=X$. The degree $0$ component of $A\otimes_KA$ is 
of dimension $1$, generated by $1\otimes 1$. An element $\lambda\cdot (1\otimes 1)$ maps 
to $1$ under the multiplication if and only if $\lambda=1$. However, 
$$
X\cdot(1\otimes 1)= (X\otimes 1)\neq  (1\otimes X)=(1\otimes 1)\cdot X
$$ 
Hence the extension is not separable. 
\item
Suppose that the characteristic of the field $K$ is different from $2$.
Consider the dg-extension $(A,d_A)\lra (B,d_B)$ of dg-division $K$-algebras,
where $(A,d_A)$ is acyclic, whence also $(B,d_B)$, and where $\ker(d_A)$ is 
a skew-field concentrated in degree $0$, and where $\ker(d_B)=D[T,T^{-1}]$
for some $T$ in non zero even degree and a skew field $D$. Then this dg-extension is not 
dg-separable. This follows from the fact that $\ker(d_B)$ is of infinite dimension over 
$\ker(d_A)$, by degree considerations, and \cite[Lemma 2.1]{ElFadilQuebec} 
shows that graded-separable extensions are finite dimensional.
Hence, the graded-extension $\ker(d_A)\lra \ker(d_B)$
is not graded-separable, and then Proposition~\ref{gradedseparableofcyclesimpliesseparable}
gives the result.
\end{itemize}
\end{Example}

We summarise the results in the following

\begin{Theorem}\label{maintheorem}
Let $K$ be a field, let $(A,d_A)$ and let $(B,d_B)$ be graded commutative dg-division rings over $K$.
Let $\varphi:(B,d_B)\ra (A,d_A)$ be a dg-extension. 
\begin{itemize}
\item Then $\ker(d_A)$ and $\ker(d_B)$ are graded-commutative graded-division rings.
\item If $(B,d_B)$ is acyclic, then also $(A,d_A)$ is acyclic and 
\begin{itemize}
\item the dg-extension is dg-separable if 
the induced graded-extension $\ker(d_B)\ra\ker(d_A)$ is graded separable. 
\item If the characteristic of $K$ is different from $2$, then the dg-extension is dg-separable if and only if 
the induced graded-extension $\ker(d_B)\ra\ker(d_A)$ is graded-separable. 
\end{itemize}
\end{itemize}
Suppose now that the characteristic of $K$ is different from $2$.
\begin{itemize}
\item If $d_B=0$ and $(A,d_A)$ is acyclic, then 
the dg-extension is not dg-separable. 
\item
If $(B,d_B)$ is acyclic and $\ker(d_B)$ is concentrated in degree $0$, and if $\ker(d_A)$ 
is not concentrated in degree $0$, then the dg-extension is not dg-separable.
\item
If  $(C,d_C)$ is a graded-commutative acyclic gr-division ring such that $\ker(d_C)$ is not concentrated in degree $0$, then there is a 
field $D$ such that $C\simeq D[T,T^{-1}]$ for some $T$ in non zero even degree. 
An extension  $D_1[T^n,T^{-n}]\ra D_2[T,T^{-1}]$ is graded-separable if and only if the field extension 
$D_1\lra D_2$ is separable and $n$ is invertible in $D_1$.
\end{itemize}
\end{Theorem}

\section{General consequences of dg-separability}

\label{generalpropertiesofdgseparableext}

\begin{Rem}
Recall that we have two concepts of semisimplicity. 
An abelian category $\mathcal A$ is semisimple if every short exact sequence 
of objects in $\mathcal A$ is split. An (graded) algebra $A$ is $J$-semisimple (Jacobson-semisimple) 
if every graded $A$-module is a direct sum of simple (graded) $A$-modules. 
It is well-known that if $A$ is artinian, then the two concepts coincide for $\mathcal A$
being the category of finitely generated (graded) $A$-modules.
Similar concepts hold for dg-modules instead of graded modules. 
\end{Rem}

\begin{Rem}\label{everyobjectprojissemisimple}
Let $\mathcal C$ be an abelian category in which every object is projective. Then $\mathcal C$
is semisimple in the sense that every short exact sequence of objects in $\mathcal C$
splits. 
\end{Rem}

\begin{Theorem} \label{acyclicimpliesequivofcategories}
Let $(A,d)$ be a dg-algebra. 
\begin{itemize}
\item \cite[Proposition 3.3]{Tempest-Garcia-Rochas} 
If a dg-module $(M,\delta)$ over $(A,d)$ is a projective 
object in the category of dg-modules,  then $(M,\delta)$ is acyclic. 
\item \cite[Theorem 4.7]{Tempest-Garcia-Rochas} 
If $(A,d)$ is acyclic, then every dg-module over $(A,d)$ is acyclic
and the functor $$A\otimes_{\ker(d)}-:gr-\ker(d)-mod\lra dg-(A,d)-mod$$
is an equivalence with quasi-inverse being the functor taking cycles. 
\item \cite[Definition 5.1 and Theorem 5.3]{Tempest-Garcia-Rochas} 
The category of dg-modules over $(A,d)$ is J-semisimple if and only if 
$(A,d)$ is acyclic and $\ker(d)$ is graded-J-semisimple.  
\end{itemize}
\end{Theorem}

We consider consequences which can be derived for dg-separable dg-extensions of dg-algebras.

\begin{Theorem} \label{separableimplieseveryobjectisdirectfactoroffree}
Let $(A,d_A)$ is a dg-algebras over some graded commutative acyclic
dg-division ring $(K,d_K)$ and suppose that $\varphi:(K,d)\ra(A,d_A)$ 
is a dg-separable dg-extension. Let $(L,d_L)$ be a graded commutative 
dg-division ring being a dg-extension of $(K,d_K)$.

Then, any dg-module $(M,\delta_M)$ over $(A\otimes_KL,d_{A\otimes_KL})$  
is a direct summand of  $(A\otimes_KL,d_{A\otimes_KL})^I$ for some index set $I$.
More precisely, $I$ is a $\ker(d_L)$-basis of $\ker(\delta_M)$. 
\end{Theorem}

Proof. 
$$\mu:A\otimes_KA\lra A$$
is split as dg-morphism by $\rho:A\lra A\otimes_KA$, satisfying $\mu\circ\rho=\id_A$. 
Then $\rho\otimes\id_L$ is a split of 
$$\mu_L:(A\otimes_KL)\otimes_L(A\otimes_KL)\lra (A\otimes_KL).$$
Indeed, 
$$(A\otimes_KL)\otimes_L(A\otimes_KL)=(A\otimes_KA)\otimes_KL$$
and with this identification we get
$$(\rho\otimes\id_L)\circ(\mu\otimes\id_L)=(\rho\circ\mu)\otimes\id_L=\id_A\otimes\id_L=\id_{A\otimes_KL}.$$
We therefore may assume that $K=L$ from the beginning. 

Doing so 
$$\mu\otimes_A\id:A\otimes_KM=A\otimes_KA\otimes_AM\lra A\otimes_AM=M$$
is split by  $\rho\otimes\id$. 
Hence $(M,\delta_M)$ is a direct factor of $(A\otimes_KM,\delta_{A\otimes_KM})$.  

We need to analyze $(A\otimes_KM,\delta_{A\otimes_KM})$.

%

Since $(K,d_K)$ is an acyclic dg-division algebra, since $(A,d_A)$ is a dg-module over $(K,d_K)$, 
as well as $(M,\delta)$, we get that $(A,d_A)$ and $(M,\delta_M)$ 
are acyclic (cf Theorem~\ref{acyclicimpliesequivofcategories}). 
Then 
$$
\Phi_K:K\otimes_{\ker(d_K)}-:\ker(d_K)-gr-mod\lra (K,d_K)-dg-mod
$$
is an equivalence of categories with inverse the functor given by taking cycles. 
Hence, the unit $\id\lra \Phi_K\circ\Phi_K^{-1}$ is an isomorphism of functors.
Moreover,
$$\Phi_A:A\otimes_{\ker(d_A)}-:\ker(d_A)-gr-mod\lra (A,d_A)-dg-mod$$
is an equivalence of categories. 
Then 
\begin{eqnarray*}
A\otimes_KM&\simeq&A\otimes_K(\Phi_K\circ\Phi_K^{-1}M)\\
&=&A\otimes_K (K\otimes_{\ker(d_K)}\ker(\delta_M))\\
&=&A\otimes_{\ker(d_K)}\ker(\delta_M)
\end{eqnarray*}
Since $\ker(d_K)$ is a gr-field, by \cite[Lemma 1.7]{vanOystaenGeel} $\ker(\delta_M)$ has a $\ker(d_K)$-basis $I$ of 
homogeneous elements. Hence, $A\otimes_KM=A^I$. 
\dickebox

\begin{Cor}\label{dg-separable-graded-semisimple}
Let $(A,d_A)$ be an acyclic dg-algebras over some graded commutative acyclic
dg-division ring $(K,d_K)$ and suppose that $\varphi:(K,d)\ra(A,d_A)$ 
is a dg-separable dg-extension.
Then any dg-module over $(A,d_A)$ is a projective object
in the category of dg-modules over $(A,d_A)$. 
Moreover, if $(A,d_A)$ is dg-artinian, then the category of dg-modules over $(A,d)$ is 
semisimple and $\ker(d_A)$ is graded-semisimple. 
\end{Cor}

Proof. 
Indeed, by Theorem~\ref{separableimplieseveryobjectisdirectfactoroffree} every 
dg-module over $(A,d_A)$ is a direct factor of $(A,d_A)^I$ for some index set $I$. 
By Theorem~\ref{acyclicimpliesequivofcategories}, since $(A,d_A)$ is assumed to be acyclic, 
$(A,d_A)$ is a projective object in the category of dg-modules over $(A,d_A)$.
Remark~\ref{everyobjectprojissemisimple} shows that this implies that the
category of dg-modules over $(A,d_A)$ is semisimple.
By \cite[Theorem 2.1]{dgHopkins}  a dg-artinian acyclic dg-algebra is dg-Noetherian. 
For dg-artinian 
and dg-Noetherian algebras
the concepts of semisimplicity and of J-semisimplicity coincides for finitely generated 
dg-modules. 
But by Theorem~\ref{acyclicimpliesequivofcategories} we get that this implies
that $\ker(d_A)$ is graded-semisimple.  \dickebox

\medskip

We can prove an analogue to \cite[Proposition 1.3]{Kuelshammernilpotent}. 
Recall that for $B\otimes_AB$-bimodules $M_1$ and $M_2$ we denote by 
$M_1^B$ the subset of elements $x$ in $M_1$ with $bx=xb$ for all $b\in B$, 
and likewise for $M_2$. 
Then for a homomorphism $\alpha:M_1\lra M_2$ of $B\otimes_AB$-bimodules we get that 
$\alpha(M_1^B)\subseteq M_2^B$. Indeed, 
$$b\alpha(x)=\alpha(bx)=\alpha(xb)=\alpha(x)b$$
for all $x\in M_1^B$ and $b\in B$. 

\begin{Theorem}\label{relativelysplit}
Let $(A,d_A)$ and $(B,d_B)$ be dg-algebras over some 
graded-commutative dg-division ring $(K,d_K)$ and suppose that 
$\varphi:(A,d_A)\lra(B,d_B)$ is a dg-separable dg-extension. 
Then any short exact sequence of dg-modules
$$0\lra (L,\delta_L)\stackrel{f}{\lra} (M,\delta_M)\stackrel{g}{\lra} (N,\delta_N)\lra 0$$
over $(B,d_B)$ is split if and only if it is split considered as a sequence of dg-modules over $(A,d_A)$.   
\end{Theorem}

Proof. The space $\Hom_{(K,d_K)}^\bullet((N,\delta_N), (M,\delta_M))$
is a dg-bimodule over $(B,d_B)-(B,d_B)^{op}$ given by $b_1\otimes b_2$
acts on $\Phi\in \Hom_{(K,d_K)}^\bullet((N,\delta_N), (M,\delta_M))$ by 
$((b_1\otimes b_2)\cdot\Phi)(n)=b_1\Phi(b_2n)$. 

\medskip

Suppose that the sequence is split as dg-modules over $(A,d_A)$. 
Let $\rho$ be an $(A,d_A)$-splitting
of the  epimorphism $g$ on the right. Then  
\begin{eqnarray*}
B\otimes_AB&\stackrel{\sigma}{\lra}&\Hom_{(K,d_K)}^\bullet((N,\delta_N), (M,\delta_M))\\
b_1\otimes b_2&\mapsto&b_1\rho b_2
\end{eqnarray*}
is a dg-$B\otimes_AB$-module homomorphism since $\rho$ is $A$-linear. 
Let $\omega\in \ker(d_{B\otimes_AB})$ mapping to $1$ under the multiplication $\mu$, 
such that $b\omega=\omega b$ for all $b\in B$ from Proposition~\ref{dgcasimir}.
But then we claim that 
$$\tau:=\sigma(\omega)\in \Hom_{(B,d_B)}^\bullet((N,\delta_N), (M,\delta_M)).$$
Indeed, 
$\sigma$ is a homomorphism of $B-B$-bimodules, and hence 
$$\tau\in \left(\Hom_{(K,d_K)}^\bullet((N,\delta_N), (M,\delta_M))\right)^B$$
by the remarks preceding the statement of Proposition~\ref{relativelysplit}.
But 
$$\left(\Hom_{(K,d_K)}^\bullet((N,\delta_N), (M,\delta_M))\right)^B=
\Hom_{(B,d_B)}^\bullet((N,\delta_N), (M,\delta_M)).$$
Further, for all $n\in N$ we get
$$(g\circ\tau)(n)=(g\circ\sigma(\omega))(n)=((g\circ\sigma)(\omega))(n)=\mu(\omega)\cdot n=1\cdot n=n
$$

\medskip

If the short exact sequence is split as sequence of dg-$(B,d_B)$-
modules, then trivially it is split as sequence of $(A,d_A)$-modules as well. 
This proves the statement. \dickebox

\begin{Rem}
Note that by Proposition~\ref{cyclestoacyclic} if $K$ is 
a graded commutative $\Z$-graded-division ring, and $(A,d)$ is a dg-division algebra, such that
  $(K,0)\lra (A,d)$ is a dg-separable dg-extension, then.the characteristic of $K$ is $2$
or $A$ cannot be graded commutative.  
\end{Rem}

Recall the concept of a separable functor introduced by 
N\u ast\u asescu, van den Bergh, and van Oystaeyen \cite{NastacescuvdBerghvOystaen}. 

\begin{Def} \cite{NastacescuvdBerghvOystaen} 
A {\em covariant 
functor} $F:{\mathcal C}\lra{\mathcal D}$ between categories $\mathcal C$ and $\mathcal D$ 
is called {\em separable} if the canonical map $\Phi^F_{A,B}:{\mathcal C}(A,B)\lra{\mathcal D}(FA,FB)$
is a naturally split monomorphism. 
\end{Def}

\begin{Prop}
Let $(A,d_A)$ and $(B,d_B)$ be dg-algebras and let $\varphi:(A,d_A)\ra (B,d_B)$ be a 
dg-extension. Then $\varphi$ is dg-separable if and only if 
the restriction $dg-(B,d_B)-mod\ra dg-(A,d_A)-mod$ is a separable functor. 
\end{Prop}

Proof. Suppose that the multiplication $B\otimes_AB\ra B$ is split.
Then the restriction is right adjoint to the induction $B\otimes_A-$ 
as is well-known (cf Yekutieli~\cite[12.6.5]{Yekutielibook}).
The counit of the adjoint pair is the multiplication map
$B\otimes_AB\ra B$. By \cite[2.2.(ii)]{WisbauerSeparability} this is equivalent with 
the fact that the restriction functor is separable. \dickebox

\end{document}